\documentclass[12pt,a4paper, leqno]{article}
\usepackage{amssymb, amscd, amsthm, amsmath, amstext, mathrsfs}
\usepackage[all]{xy}

\usepackage[colorlinks, linkcolor=blue, anchorcolor=green, citecolor=blue]{hyperref}

\newtheorem{thm}{Theorem}[section]
\newtheorem{coro}[thm]{Corollary}
\newtheorem{lemma}[thm]{Lemma}
\newtheorem{prop}[thm]{Proposition}

\theoremstyle{remark}
\newtheorem{remark}[thm]{\textbf{Remark}}

\theoremstyle{definition}

\numberwithin{equation}{thm}
\newcommand{\newpara}{\noindent\refstepcounter{thm}{\bf(\thethm)\;}}

\newcommand{\CH}{\mathrm{CH}}
\newcommand{\tors}{\mathrm{tors}}
\newcommand{\vp}{\varphi}

\newcommand{\dgf}[1]{\langle\, #1\,\rangle}
\newcommand{\pfi}[1]{\langle\!\langle\, #1\,\rangle\!\rangle}

\newcommand{\sK}{\mathscr{K}}

\def\ker{\mathrm{Ker}}

\newcommand{\N}{\mathbb N}
\newcommand{\Q}{\mathbb{Q}}
\newcommand{\Z}{\mathbb Z}

\newcommand{\al}{\alpha}

\newcommand{\car}{\mathrm{char}}

\newcommand{\ind}{\mathrm{ind}}
\newcommand{\Br}{\mathrm{Br}}
\newcommand{\lra}{\longrightarrow}
\newcommand{\ov}[1]{\overline{#1}}
\newcommand{\simto}{\xrightarrow{\sim}}

\newcommand{\Pfi}[1]{\langle\!\langle\, #1\,]\!]}

\newcommand{\A}{\mathbb{A}}
\newcommand{\bH}{\mathbb{H}}

\begin{document}
\title{\textbf{Chow Groups of Quadrics in Characteristic Two}}
\author{Yong Hu, Ahmed Laghribi, and Peng Sun}

\maketitle

\begin{abstract}
Let $X$ be a smooth projective quadric defined over a field  of characteristic 2. We prove that in the Chow group of codimension 2 or 3 of $X$ the torsion subgroup has at most two elements. In codimension 2, we determine precisely when this torsion subgroup is nontrivial. In codimension 3, we show that there is no torsion if $\dim X\ge 11$. This extends the analogous results in characteristic different from 2, obtained by Karpenko in the nineteen-nineties.
\end{abstract}


\medskip

\noindent {\bf MSC classification 2020:} 11E04, 14C35, 14C25, 19E08

\section{Introduction}

Let  $X$ be a smooth projective variety over a field $k$.  For each natural number $p$, denote by $\CH^p(X)$ the Chow group of codimension $p$ cycles on $X$ modulo rational equivalence (\cite{Ful}). When $p\ge 2$, determining the structure of the group $\CH^p(X)$, especially that of the torsion subgroup, is an interesting but often difficult problem in algebraic geometry. A closely related notion is  the Grothendieck ring $K_0(X)$ of vector bundles on $X$. A consequence of the Grothendieck--Riemann--Roch theorem (cf. \cite[\S\;15.2]{Ful}) is that the Chern character gives an isomorphism $K_0(X)\otimes\Q\cong \bigoplus_{p\ge 0}\CH^p(X)\otimes\Q$.

Consider the special case where $X$ is a smooth projective quadric. Chow groups and $K$-theory of $X$ were first studied by Swan in \cite{Swan85AnnMath} and \cite{Swan89PAMS}). In the 1990's, Karpenko made a systematic study of the structure of $\CH^p(X)$ for $p\le 4$ based on Swan's work (\cite{Karpenko90AlgGeoInv}, \cite{Karpenko91ChowGroups}, \cite{Karpenko91cyclesCodim3}, \cite{Karpenko95Nova}, \cite{Karpenko96DocMath}; see also \cite{KarpenkoMerkurjev90ChowGroupsAlgiAna}). While Swan's main theorem (\cite[Thm.\;1]{Swan85AnnMath}) on the $K$-theory of $X$ is established in arbitrary characteristic, Karpenko's theorems are stated only in characteristic different from 2. Among others he proves the following results in codimensions 2 and 3 (in characteristic $\neq 2$):

\begin{enumerate}
  \item (\cite[Thm.\;6.1]{Karpenko90AlgGeoInv}) The torsion subgroup $\CH^2(X)_{\tors}$ of $\CH^2(X)$ is finite of order at most 2, and it is nontrivial if and only if the quadratic form defining $X$ is an anisotropic 3-fold Pfister neighbor.

      In particular, $\CH^2(X)_{\tors}=0$ if $\dim X>6$.
  \item The torsion subgroup $\CH^3(X)_{\tors}$ of $\CH^3(X)$ is finite of order at most 2 (\cite{Karpenko91cyclesCodim3}). It is trivial if $\dim X>10$ (\cite[Thm.\;6.1]{Karpenko95Nova}). The proof of the latter result depends on a theorem of Rost about 14-dimensional forms with trivial discriminant and trivial Clifford invariant (\cite{Rost99dim14}, \cite{Rost06}). Without using Rost's theorem, one can show $\CH^3(X)_{\tors}=0$ when $\dim X>14$ (\cite[Thm.\;7.1]{Karpenko95Nova}).
\end{enumerate}

Still in characteristic different from $2$,  Izhboldin has further developed Karpenko's methods and obtained more precise information about $\CH^3(X)_{\tors}$ when $7\le \dim X\le 10$ (cf. \cite[Thm.\;0.5]{Izhboldin01}).

As applications, Karpenko's results on Chow groups have been used by Kahn, Rost and Sujatha  to compute the unramified cohomology groups up to degree 4 for projective quadrics (\cite{Kahn95ArchMath}, \cite{KahnRostSujatha98}, \cite{KahnSujatha00JEMS}, \cite{KahnSujatha01Duke}). In turn,  some of their computations make it possible to prove similar results for codimension 4 Chow groups (\cite[\S\;8]{Karpenko95Nova}, \cite{Karpenko96DocMath}) and some others play a key role in Izhboldin's construction of a field of $u$-invariant 9 in characteristic $\neq 2$ (\cite[Thm.\;0.1]{Izhboldin01}).

\

It is natural to expect the same results as above in characteristic 2.  Basically, one can follow the same methods as in Karpenko's papers. But on the one hand, at some points the original proofs need appropriate modifications, where quite a few details are worth clarifying with special care. On the other hand, it does seem (at least to us) that some other arguments in Karpenko's work (for example, those in \cite[\S\;6]{Karpenko90AlgGeoInv} and \cite[\S\;6]{Karpenko95Nova})  rely on quadratic form techniques  which are particular in characteristic different from 2. In their construction of a field with a special indecomposability phenomenon, Barry, Chapman and the second named author have used the vanishing of $\CH^2(X)_{\tors}$ in characteristic 2 (\cite[Thm.\;A.1]{BarryChapmanLaghribi20Israel}). Their proof of this vanishing result provides an example of adapting Karpenko's arguments in characteristic 2.

\medskip

In this paper, we make a further study of the Chow groups $\CH^2(X)$ and $\CH^3(X)$ in characteristic 2 and extend the other results of Karpenko mentioned above. In particular, we show that the group $\CH^3(X)_{\tors}$ has at most two elements, as in the case of characteristic $\neq 2$. {We also prove $\CH^3(X)_{\tors}=0$ as soon as $\dim X>10$ (Theorem\;\ref{7p11}). Here we need to extend Rost's theorem to characteristic 2 (Theorem\;\ref{7p10}), which we do using a specialization argument, and which we believe has independent interest.

Karpenko's results for $\CH^4(X)_{\tors}$ and {some of} Izhboldin's results for $\CH^3(X)_{\tors}$ rely on the computation of degree 4 unramified cohomology groups established in \cite{KahnRostSujatha98}. As we would like to leave out discussions on unramified cohomology in this paper, {we will not investigate full generalizations of these results in characteristic 2. We only provide a few examples where $\CH^3(X)_{\tors}=0$ for some quadrics in lower dimensions (see Prop.\;\ref{7p12} and Remark\;\ref{7p16}).} A study of $\CH^4(X)_{\tors}$ is likely to be the theme of a further work.

It is perhaps worth mentioning that to obtain the main result (Theorem\;\ref{5p3}) for $\CH^2(X)_{\tors}$, we give a proof that works uniformly for quadrics of all dimensions $\ge 3$. Unlike the methods in \cite{Karpenko90AlgGeoInv} and the appendix of \cite{BarryChapmanLaghribi20Israel} (see Remark\;\ref{5p4} and the appendix), this approach utilizes Kato--Milne cohomology (and also $K$-cohomology implicitly).

In \cite{HuSun23}, the results of this paper have been used to study unramified cohomology in degree $\le 3$ for quadrics in characteristic $2$.

\

\noindent {\bf Notation.} Throughout what follows, $k$ denotes a field of characteristic 2, with a fixed separable closure $\ov{k}$. Let $\wp(k)$ be the image of the Artin--Schreier map $\wp: k\to k\,;\; x\mapsto x^2-x$.

For an algebraic $k$-variety $Y$,  we write $\ov{Y}=Y\times_k\ov{k}$.

For an abelian group $M$, we denote by $M_{\tors}$ its torsion subgroup.

\section{Preliminaries on Chow groups of quadrics}

\newpara\label{2p1} We assume the reader is familiar with the general theory of quadratic forms.
Terminology and notation about quadratic forms, if not explained, are standard and are adopted from \cite{EKM08}. In particular, a quadratic form is called \emph{nondegenerate} if the radical of its polar bilinear form is anisotropic of dimensional $\le 1$. A quadratic form of dimension $\ge 2$ is nondegenerate if and only if its projective quadric  is smooth as an algebraic variety.

A nondegenerate quadratic form $\vp$ over $k$ has the following normal form
\[
\begin{split}
\vp&\cong [a_1,\,b_1]\bot \cdots\bot [a_m,\,b_m]\,,\quad \text{ if }\dim\vp=2m\,,\\
\vp&\cong [a_1,\,b_1]\bot \cdots\bot [a_m,\,b_m]\bot \langle c\rangle\,,\quad \text{ if } \dim\vp=2m+1\,,
\end{split}
\]
where $a_i,\,b_i\in k$, $c\in k^*$, and $[a_i, b_i]$ (resp. $\langle c\rangle$) denotes the quadratic form $a_ix^2+xy+b_iy^2$ (resp. $cx^2$). In the even dimensional case, the \emph{Arf invariant} (or \emph{discriminant}) of $\vp$ is defined as the image of the element $\sum^m_{i=1}a_ib_i$ in the quotient group $k/\wp(k)$. It is uniquely determined by $\vp$ and denoted by $\mathrm{Arf}(\vp)$.  The $k$-algebra $k[T]/(T^2-T-\al)$, where $\al\in k$ is a representative of the Arf invariant $\mathrm{Arf}(\vp)\in k/\wp(k)$, is uniquely determined. It will be called the \emph{discriminant algebra} of $\vp$.

Let $W_q(k)$ be the Witt group of even-dimensional nondegenerate quadratic forms over $k$.  For any integer $n\geq 1$, let $I^n_q(k)$ be the subgroup of $W_q(k)$ generated by $n$-fold quadratic Pfister forms.

We will need the Arason-Pfister Hauptsatz, simply called the Hauptsatz, that asserts the following: If an anisotropic quadratic form $\vp$ belongs to $I^n_q(k)$, then it has dimension $\geq 2^n$ (\cite[Thm.\;23.7]{EKM08}, \cite[Thm.\;4.2 (iv)]{HoffmannLaghribi04TAMS}).

For two quadratic forms $\vp$ and $\psi$ over $k$, we say $\psi$ is a \emph{subform} of $\vp$ if $\psi\cong\vp|_{W}$ for some subspace $W$ in the vector space $V_{\vp}$ of $\vp$.
{When this happens we write $\psi\subseteq\vp$.}

\medskip

\newpara\label{2p2} Now we recall some known facts about Chow groups of projective quadrics (which are valid in arbitrary characteristic). More details can be found in \cite[\S\;2]{Karpenko90AlgGeoInv} and \cite[\S\;68]{EKM08}.

Let $\vp$ be a nondegenerate quadratic form of dimension $\ge 3$ over $k$, defined on a $k$-vector space $V$. Let $X=X_{\vp}$ be the projective quadric defined by $\vp$, which is a closed subvariety in the projective space $\mathbb{P}(V)$. Let $h\in\CH^1(X)$ be the pullback of the class of a hyperplane in $\mathbb{P}(V)$. For each $p\in\N$, the power $h^p$ generates a torsion free subgroup $\Z.h^p$ in $\CH^p(X)$, called the \emph{elementary part} of $\CH^p(X)$. We say $\CH^p(X)$ is \emph{elementary} if it is equal to its elementary part.

\

The following result is well known.

\begin{prop}\label{2p3}
  Let $\vp$ be a nondegenerate quadratic form of dimension $\ge 3$ over $k$, and let $X_{\vp}$ be the projective quadric defined by $\vp$. Assume that $\vp$ is anisotropic.

  \begin{enumerate}
    \item If $\dim\vp\neq 2p+2$, then $\CH^p(X_{\vp})$ is elementary if and and only if $\CH^p(X_{\vp})$ is torsion free.
    \item If $\dim\vp=2p+2$, then $\CH^p(X_{\vp})$ is elementary if and and only if $\CH^p(X_{\vp})$ is torsion free and $\mathrm{Arf}(\vp)\neq 0$.
  \end{enumerate}
\end{prop}

We have some known examples of torsion free Chow groups.

\begin{prop}\label{2p4}
  Let $X$ be a smooth projective quadric of dimension $d\ge 1$ over $k$.

  \begin{enumerate}
    \item The groups $\CH^0(X),\,\CH^1(X)$ and $\CH^d(X)$ are torsion free.
    \item If $X$ is isotropic, then $\CH^2(X)$ is torsion free.
  \end{enumerate}
\end{prop}
\begin{proof}
The statement for $\CH^d(X)$ follows from \cite[Cor.\;70.4]{EKM08}. (See also \cite[Lemma\;4.1]{Totaro08JAG} for a generalization in the singular case. In characteristic $\neq 2$ this was proved independently in \cite{Swan89PAMS} and \cite[Prop.\;2.6]{Karpenko90AlgGeoInv}.)
The other statements are well known.
\end{proof}

\section{Clifford algebra and splitting index}

Throughout this section, let $\vp$ be a nondenegerate quadratic form of dimension $\ge 1$ over $k$, and let $C(\vp)$ and $C_0(\vp)$ be its Clifford algebra and even Clifford algebra respectively.

\medskip

\newpara\label{3p1} In the sequel we will frequently use a simple $k$-algebra $C_0'(\vp)$ defined as follows: If $\vp$ has even dimension and trivial Arf invariant, then $C_0(\vp)\cong A\times A$ for a unique (up to isomorphism) central simple $k$-algebra $A$ and $C(\vp)\cong \mathrm{M}_2(A)$ (cf. \cite[Remark\;13.9]{EKM08}). In this case we set $C_0'(\vp)=A$. Otherwise ($\dim\vp$ is odd, or $\dim\vp$ is even but $\mathrm{Arf}(\vp)\neq 0$), we put  $C'_0(\vp)=C_0(\vp)$.

In any case, we can write $C'_0(\vp)\cong \mathrm{M}_{2^s}(D)$ for some $s\in \N$ and some division algebra $D$ with the same center as $C'_0(\vp)$. We write $s(\vp)$ for the integer $s$ here and define $\ind(\vp)=\ind(C'_0(\vp))$, the Schur index of $C'_0(\vp)$ over its center. Following \cite{Karpenko90AlgGeoInv} and \cite{Izhboldin01}, we call $\ind(\vp)$ and $s(\vp)$ the \emph{index} and the \emph{splitting index} of $\vp$ respectively.

From the definitions we find easily the relation
\begin{equation}\label{3p1p1}
  s(\vp)+\log_2\bigl(\ind(\vp)\bigr)=\left[\frac{\dim\vp-1}{2}\right]\;.
\end{equation}Also, it is easy to see
\begin{equation}\label{3p1p2}
  \begin{cases}
 i_W(\vp)\le s(\vp)\le \left[\frac{\dim\vp-1}{2}\right] \quad & \text{ if } \vp \text{ is not hyperbolic}\,,\\
 s(\vp)=i_W(\vp)-1=\left[\frac{\dim\vp-1}{2}\right] \quad & \text{ if } \vp \text{ is hyperbolic}\,,\\
\end{cases}
\end{equation}
where $i_W(\vp)$ denotes the Witt index of $\vp$.

\

We have some auxiliary results where the splitting index is used to detect the structure of quadratic forms in low dimensions.

\begin{lemma}\label{3p2}
  Suppose $\dim\vp=5$ (so that $0\le s(\vp)\le 2$ by $\eqref{3p1p2}$).

  \begin{enumerate}
    \item  $s(\vp)=2\iff i_W(\vp)=2$.
    \item  Assume that $\vp$ is anisotropic. Then the following conditions are equivalent:

    \begin{enumerate}
      \item $s(\vp)=1$.

      \item For some quadratic separable extension $K/k$, the form $\vp_K$ splits completely, i.e., $i_W(\vp_K)=2$.
      \item $\vp$ is similar to $\psi\bot\langle c\rangle$ for some $c\in k^*$ and some $2$-fold Pfister form $\psi$.

      \item $\vp$ is a Pfister neighbor.
    \end{enumerate}
  \end{enumerate}
\end{lemma}
\begin{proof}
Let us write $\vp=\psi\bot \langle c\rangle$ with $c\in k^*$. Then we have $C(-c\psi)=C(c\psi)\cong C_0(\vp)$ (noticing that $\car(k)=2$) and hence $\ind(\vp)=\ind(C(c\psi))$.  From \eqref{3p1p1} we see
\[
 \begin{split}
 s(\vp)=0&\iff C(c\psi) \text{ is a central division $k$-algebra of degree } 4\,,\\
   s(\vp)=1&\iff C(c\psi)\cong \mathrm{M}_2(Q) \text{ for some quaternion division $k$-algebra } Q\,,\\
    s(\vp)=2&\iff C(c\psi)\cong \mathrm{M}_4(k)=C(2\mathbb{H})\;.
 \end{split}
 \]

If $i_W(\vp)=2$, then clearly \eqref{3p1p2} yields $s(\vp)=2$. Conversely, if $s(\vp)=2$, then  we have $C_0(c\psi)\cong C_0(2\mathbb{H})$. Hence, by \cite[\S\,9, Thm.\;7]{Knus88book}, $c\psi$ is similar to $2\mathbb{H}$, i.e., $\psi$ is hyperbolic, giving $i_W(\vp)=2$. This proves (1).

\medskip

To prove (2), let us consider the Albert form $\gamma:=\psi \perp c[1, r]$, where $r$ is a representative of the Arf invariant of $\psi$. By  Jacobson's theorem \cite{MS} one has $s(\varphi)=1$ if and only if $\gamma$ is isotropic. Moreover, if $\gamma$ is isotropic, then $\gamma \cong \tau \perp \mathbb{H}$ for some form $\tau$ similar to a $2$-fold Pfister form. Adding the form $\langle c\rangle$ in both sides yields $\tau \perp \langle c\rangle \perp \mathbb{H} \cong \psi \perp \langle c\rangle \perp \mathbb{H}$. Cancelling the hyperbolic plane, we get $\varphi \cong \tau \perp \langle c\rangle$, this proves (a)$\Rightarrow$(c). The implication (c)$\Rightarrow$(a) is clear from the definition of $s(\varphi)$. The equivalence (c)$\Leftrightarrow$(d) is stated in \cite[Prop.\;3.2\,(3)]{L}. Since any $2$-fold Pfister form is split by a separable quadratic extension, the implication (c)$\Rightarrow$(b) is immediate. For the implication (b)$\Rightarrow$(c) we use the fact that if an anisotropic quadratic form becomes isotropic over a separable quadratic extension $K=F[X]/(X^2-X-a)$, then the form contains a subform similar to $[1, a]$ (\cite[Prop.\;34.8]{EKM08}).
\end{proof}

\begin{lemma}\label{3p3}
  Suppose $\dim\vp=6$.  (Thus $0\le s(\vp)\le 2$ by $\eqref{3p1p2}$.)

  \begin{enumerate}
    \item Assume that $\vp$ is an Albert form, i.e.,  $\mathrm{Arf}(\vp)=0$. Then
    \[
    \begin{split}
      s(\vp)=0&\iff \ind(C(\vp))=4 \iff i_W(\vp)=0\,,\\
       s(\vp)=1&\iff \ind(C(\vp))=2 \iff i_W(\vp)=1\,,\\
        s(\vp)=2&\iff \ind(C(\vp))=1 \iff i_W(\vp)=3\,.\\
    \end{split}
    \]
    \item Assume that $\mathrm{Arf}(\vp)\neq 0$, so that the discriminant algebra of $\vp$ is a separable quadratic field extension $K$ of $k$.

    Then
     \[
    \begin{split}
      s(\vp)=0&\iff \ind(C_0(\vp))=4 \iff i_W(\vp_K)=0\,,\\
       s(\vp)=1&\iff \ind(C_0(\vp))=2 \iff i_W(\vp_K)=1\,,\\
        s(\vp)=2&\iff \ind(C_0(\vp))=1 \iff i_W(\vp_K)=3\,.\\
    \end{split}
    \]

  \end{enumerate}
\end{lemma}
\begin{proof}
    Combine \eqref{3p1p2} with   \cite[\S\,11,\, Cor.\;5 and Remark\;13]{Knus88book}  (see also \cite[(16.5)]{KMRT}).
\end{proof}

The following lemma includes a characteristic 2 version of \cite[(5.4)]{Karpenko90AlgGeoInv}. It can be proved in a similar way, with the help of Lemma\;\ref{3p3} and \cite[Prop.\;34.8]{EKM08}.

\begin{lemma}\label{3p4}
  Suppose $\dim\vp=6$. Let $Z$ be the discriminant algebra of $\vp$ and let $N_{Z/k}: Z\to k$ denote the norm of $Z/k$ regarded as a binary quadratic form.

  Assume that $\vp$ is anisotropic.

    \begin{enumerate}
    \item  The following conditions are equivalent:

   \begin{enumerate}
   \item $\vp\cong \dgf{a,\,b,\,c}_{\mathrm{bil}}\otimes N_{K/k}$ for some $a,\,b,\,c\in k^*$ and some quadratic separable extension $K/k$.

   Here $N_{K/k}: K\to k$ denotes the norm considered as a binary quadratic form, and $\dgf{a,\,b,\,c}_{\mathrm{bil}}$ denotes the ternary bilinear form $ax_1x_2+by_1y_2+cz_1z_2$.
     \item $\vp$ is similar to $\dgf{1,\,a,\,b}_{\mathrm{bil}}\otimes N_{K/k}$ for some $a,\,b\in k^*$ and some quadratic separable extension $K/k$.
     \item  $\vp$ is a Pfister neighbor.
     \item $\vp_K$ is hyperbolic for some quadratic separable extension $K/k$.
     \item $\vp$ is not an Albert form and $s(\vp)=2$.
   \end{enumerate}
  Note that when the above conditions hold, $K/k$ must be the discriminant algebra of $\vp$ and $\vp$ has a decomposition $\vp=\psi\bot\theta$, where $\psi=\dgf{a\,,\,b}_{\mathrm{bil}}\otimes N_{K/k}$ is similar to a $2$-fold Pfister form.
     \item Suppose $\mathrm{Arf}(\vp)\neq 0$. Then the following are equivalent:

     \begin{enumerate}
       \item $s(\vp)=1$ (i.e., $C_0(\vp)\cong \mathrm{M}_2(Q)$ for some quaternion division $Z$-algebra $Q$).
       \item $i_W(\vp_Z)=1$.
       \item $\vp\cong c.N_{Z/k}\bot\psi$,  where $c\in k^*$, $\psi$ is similar to a $2$-fold Pfister form and $\psi_Z$ is anisotropic.
     \end{enumerate}

      \item Suppose $\mathrm{Arf}(\vp)\neq 0$. Then the following are equivalent:

     \begin{enumerate}
       \item $s(\vp)=0$ (i.e., $C_0(\vp)$ is a central division algebra of degree $4$ over $Z$).
       \item $\vp_Z$ is anisotropic, i.e., $i_W(\vp_Z)=0$.
       \item $\vp$ cannot be written as $\psi\bot\theta$, where $\psi$ is similar to a $2$-fold Pfister form.
     \end{enumerate}

    \end{enumerate}
\end{lemma}

The first assertion in the lemma below is a characteristic 2 analogue of \cite[(3.3)]{Karpenko91cyclesCodim3}. It has been proved in \cite[Lemma\;3.6]{Laghribi15}. Taking \eqref{3p1p1} and \eqref{3p1p2} into account, we can deduce the second assertion easily.

\begin{lemma}\label{3p5}
Suppose that $\vp$ is anisotropic of dimension $8$ and $\mathrm{Arf}(\vp)=0$. (Note that $0\le s(\vp)\le 3$ by $\eqref{3p1p2}$.)

  \begin{enumerate}

    \item If $s(\vp)\ge 2$, then there exist $a,\,b,\,c\in k^*$ and a separable quadratic extension $L/k$ such that $\vp$ is similar to $(\langle 1,\,a\rangle_{\mathrm{bil}}\bot c.\langle 1,\,b\rangle_{\mathrm{bil}})\otimes N_{L/k}$. In particular, $\vp$ has a decomposition $\vp=\vp_1\bot\vp_2$ where both $\vp_1$ and $\vp_2$ are scalar multiples of $2$-fold Pfister forms.
   \item $s(\vp)=3$ if and only if $\vp$ is similar to a  $3$-fold Pfister form.
  \end{enumerate}
\end{lemma}

\section{Some $K$-theory of quadrics}

We briefly review some useful results from the $K$-theory of smooth projective quadrics, which  hold in arbitrary characteristic.

\medskip

Throughout this section, let $\vp$ be a nondegenerate quadratic form of dimension $\ge 3$ over the field $k$ and let $X=X_{\vp}$ be the smooth projective quadric defined by $\vp$.

\

\newpara\label{4p1} Let $K_0(X)$ be the Grothendieck group of isomorphism classes of coherent sheaves on $X$ modulo an equivalence relation defined via short exact sequences. The natural topological filtration on $K_0(X)$ will be denoted by $K_0(X)^{(p)},\,p\in\N$. For each $i\in\N$, we put
\[
K_0(X)^{(i/i+1)}:=\frac{K_0(X)^{(i)}}{K_0(X)^{(i+1)}}\;.
\]
By \cite{SGA6} (see also \cite[\S\,15.1]{Ful}), there is  a natural surjection
\[
  \rho^i\,:\;\CH^i(X)\lra K_0(X)^{(i/i+1)}\;;\quad [Z]\longmapsto [\mathscr{O}_Z]\;,
\]and the kernel of $\rho^i$ is torsion.  In fact, $\rho^i$  is an isomorphism for $i\in \{0,\,1,\,2,\,3,\,\dim X\}$ (\cite[(3.1)]{Karpenko90AlgGeoInv}).
For the injectivity of $\rho^3$, we can use the Brown--Gersten--Quillen spectral sequence in higher $K$-theory (cf. \cite{QuillenHighK}) and follow the ideas in the proof of \cite[(4.5)]{Karpenko90AlgGeoInv}. It is sufficient to notice that
for any field extension $E/k$, the natural map $H^1(X,\,\sK_2)\to H^1(X_E,\,\sK_2)$ is injective by \cite[Prop.\;1.5]{Merkurjev95ProcSympos58}.

 By abuse of notation, let $h$ also denote the class of the structural sheaf of a hyperplane section in $X$.
For each $i\in\N$, we say that $K_0(X)^{(i/i+1)}$ is \emph{elementary} if it is generated by the image of $h^i$.

\medskip

As already observed in \cite[Appendix]{BarryChapmanLaghribi20Israel}, we have the following variant of \cite[(4.5)]{Karpenko95Nova}:

\begin{prop}\label{4p2}
 Suppose that $\vp$ can be written as $\vp:=a[1,\,d]\bot\rho$, where $a\in k^*$ and $d\in k$ represents $\mathrm{Arf}(\rho)$. Let $\psi=\langle a\rangle\bot\rho$.

 For any $p\in\N$, if the groups $K_0(X_{\vp})^{(i/i+1)}$ are elementary for all $i\le p$, then the groups $K_0(X_{\psi})^{(i/i+1)}$ are also elementary for all $i\le p$.
\end{prop}

We also need the characteristic 2 version of \cite[(4.7)]{Karpenko95Nova}.

\begin{prop}\label{4p3}
Suppose $\vp=\rho\bot a.[1,\,b]$, where $a\in k^*$ and $\rho$ is an even-dimensional form.
  Let $\psi=\rho\bot \langle a\rangle$.  Assume that the discriminant algebra of $\vp$ is a quadratic field extension $K/k$ such that $\ind\big(C_0(\psi)_K\big)=\ind\big(C_0(\psi)\big)$.

  If for some $p\in\N$ the groups $K_0(X_{\psi})^{(i/i+1)}$ are elementary for all $i\le p-1$, then the groups $K_0(X_{\vp})^{(i/i+1)}$ are elementary for all $i\le p$.
\end{prop}

\section{Codimension two cycles on projective quadrics}

In this section we prove our main results about codimension two Chow groups.

As in the previous section, let $X=X_{\vp}$ be a smooth projective quadric of dimension $d\ge 1$, defined by a nondegenerate quadratic form $\vp$ over $k$. We will write
\[
 \CH^*(X):=\bigoplus_{i\ge 0}\CH^i(X)\quad\text{and}\quad \mathrm{Gr}K_0(X):=\bigoplus_{i\ge 0}K_0(X)^{(i/i+1)}\,.
\]By \cite[Cor.\;70.4]{EKM08}, we have $\CH^d(X)=\Z.[x]$, where $x\in X$ is  a closed point of minimal degree. (In characteristic $\neq 2$ this fact was proved independently in \cite{Swan89PAMS} and \cite[Prop.\;2.6]{Karpenko90AlgGeoInv}.) To study the Chow group $\CH^2(X)$ we need only to consider the case $d=\dim X\ge 3$.

\

First we observe that the cases with $\dim X=3$ or 4 can be treated in the same way as in \cite{Karpenko90AlgGeoInv}, using the isomorphism $\CH^*(X)\cong \mathrm{Gr}K_0(X)$ (cf. \eqref{4p1}).

\begin{thm}[{\cite[(5.3)]{Karpenko90AlgGeoInv}}]\label{5p1}
  Assume that $\vp$ is an anisotropic form of dimension $5$.

 Then  $\CH^2(X)_{\tors}\cong (\Z/2)^{s(\vp)}$ and $s(\vp)=0$ or $1$.

 Moreover, $s(\vp)=1$ if and only if $\vp$ contains a scalar multiple of  a $2$-fold Pfister form, if and only if $\vp$ is a Pfister neighbor (cf. Lemma\;$\ref{3p2}$).
\end{thm}

\begin{thm}[{\cite[(5.5)]{Karpenko90AlgGeoInv}}]\label{5p2}
  Assume that $\vp$ is an anisotropic form of dimension $6$.

  \begin{enumerate}
    \item If $\vp$ is an Albert form, i.e., $\mathrm{Arf}(\vp)=0$, then the group $\CH^*(X)=\bigoplus_{i\ge 0}\CH^i(X)$ is torsion free and $\CH^2(X)$ can be identified with the subgroup
   $\Z.h^2\oplus \Z.4\ell_2$ of $\CH^2(\ov{X})=\Z.h^2\oplus\Z.\ell_2$.

   (Here $\ell_2$ denotes the class of a $2$-dimensional linear space.)
    \item Assume that $\mathrm{Arf}(\vp)\neq 0$.

    \begin{enumerate}
      \item If $\vp$ is a Pfister neighbor, i.e., $s(\vp)=2$ (cf. Lemma$\;\ref{3p4}\;(1)$), then $\CH^3(X)_{\tors}$ and $\CH^2(X)_{\tors}$ are both isomorphic to $\Z/2$.
      \item If $s(\vp)=1$ (cf. Lemma$\;\ref{3p4}\;(2)$), then $\CH^3(X)_{\tors}\cong \Z/2$ and $\CH^2(X)_{\tors}=0$.
      \item If $s(\vp)=0$  (cf. Lemma$\;\ref{3p4}\;(3)$), then  $\CH(X)$ is torsion free.
    \end{enumerate}
  \end{enumerate}
\end{thm}

Our goal now is to prove the following:

\begin{thm}[{See \cite[(6.1)]{Karpenko90AlgGeoInv} in characteristic $\neq 2$}]\label{5p3}
 Let $X=X_{\vp}$ be the projective quadric defined by a nondegenerate quadratic form $\vp$ of dimension $\ge 5$ over $k$.

Then $\CH^2(X)_{\tors}$ is either $0$ or isomorphic to $\Z/2\Z$.

Moreover,  $\CH^2(X)_{\tors}\cong \Z/2$ if and only if $\vp$ is an anisotropic  $3$-fold Pfister neighbor.
\end{thm}

\begin{remark}\label{5p4}
  If $\car(k)\neq 2$, Karpenko's proof of Theorem\;\ref{5p3} is based on the  following observation (cf. \cite[(6.2)--(6.3)]{Karpenko90AlgGeoInv}):

\medskip

 \emph{Assume $\vp$ is nondegenerate of dimension $\ge 7$. Then there exists a purely transcendental extension $L/k$ and a nondegenerate $6$-dimensional quadratic form $\psi$ over $L$ such that the following properties hold:}

  \begin{enumerate}
    \item \emph{The transcendence degree $\mathrm{trdeg}(L/k)$ is equal to $\dim\vp-6$.}
    \item \emph{If $\vp$ is anisotropic over $k$, then $\psi$ is anisotropic over $L$.}
    \item \emph{Letting $X_{\vp}/k$ and $X_{\psi}/L$ be the projective quadrics defined by $\vp$ and $\psi$ respectively, we have} $\CH^2\left(X_{\vp}\right)_{\tors}\cong \CH^2\left(X_{\psi}\right)_{\tors}$.
    \item If $\vp$ is anisotropic, then \emph{$\psi$ is a $3$-fold Pfister neighbor if and only if  $\vp$ is a $3$-fold Pfister neighbor.}
  \end{enumerate}

Here properties (1) and (2) are clear from the construction. Karpenko verified property (3) by using excision and fibration arguments, and he proved property (4) with the help of some algebraic theory of quadratic forms in characteristic $\neq 2$. In \cite[Appendix]{BarryChapmanLaghribi20Israel}, Barry, Chapman and Laghribi have shown that Karpenko's method can be adapted to deal with the case $\dim\vp>8$ in characteristic 2. In their construction a form $\psi$ of dimension 9 (instead of 6) is used, and hence there is no need to check a condition similar to property (4) above.

\medskip

 When $\dim\vp$ is 7 or 8, we still have an adapted version of Karpenko's arguments in characteristic 2, thus obtaining a proof of Theorem\;\ref{5p3}. More precisely, following some ideas in \cite[Appendix]{BarryChapmanLaghribi20Israel}, with some extra effort we can  prove  the following in characteristic 2:

\emph{Write $\vp=\vp_0\bot [b,\,c]\bot\tau$ with $b,\,c\in k$ and $\dim\tau=4$. Then there exists a purely transcendental extension $L/k$ and an element $f\in L^*$ such that the above properties (1)--(4) hold for $\vp$ and $\psi:=[f,\,c]\bot\tau$.}

\medskip

If $\dim\vp=8$, i.e., $\vp_0=[a_0,\,a_1]$ for some $a_0,\,a_1\in k$, then similar to Claim 2 of Case 1 in \cite[Appendix, p.318]{BarryChapmanLaghribi20Israel}, we have
$\CH^2(X_{\vp})_{\tors}\cong \CH^2(X_{\vp'})_{\tors}$, where $\vp'$ is the quadratic form $\langle a_0t^2+t+a_1\rangle\bot [b,\,c]\bot\tau$ defined over $k(t)$. This allows us to reduce the 8-dimensional case to the 7-dimensional one.

Now let us assume $\dim\vp=7$. Then $\vp=\langle a\rangle\bot[b,\,c]\bot\tau$ for some $a\in k^*$. Replacing $c$ by $a+c$ if necessary, we may assume $c\neq 0$. Let $U\subseteq \mathbb{A}^6$ be the affine quadric defined by $a+bx^2+xy+cy^2+\tau=0$ and consider the projection onto the $x$-coordinate $\pi: U\to \mathbb{A}^1$. The generic fiber $U'$ of $\pi$ is the affine quadratic defined by $f+y+cy^2+\tau=0$ over the rational function field $L:=k(x)$, where $f:=(a+bx^2)x^{-2}\in L=k(x)$. Using excision sequences we see that $\CH^2(U)\cong \CH^2\left(X_{\vp}\right)_{\tors}$ and $\CH^2(U')\cong \CH^2\left(X_{\psi}\right)_{\tors}$. To show the desired isomorphism $\CH^2\left(X_{\vp}\right)_{\tors}\cong \CH^2\left(X_{\psi}\right)_{\tors}$, it suffices to prove that for every closed point $P\in\A^1$, the fiber $U_P$ of $\pi$ over $P$ satisfies $\CH^1(U_P)=0$. If $x(P)=0$, this follows from Claim 1 of Case 2 in \cite[Appendix, p.318]{BarryChapmanLaghribi20Israel}. If $x(P)\neq 0$, we can use Lemma\;\ref{6p3}.

We have thus constructed the form $\psi=[f,\,c]\bot \tau$ such that properties (1), (2) and (3) hold. We prove property (4) in the appendix.
\end{remark}

Our goal now is to give a proof of Theorem\;\ref{5p3} which works in a uniform way for all nondegenerate forms of dimension $\ge 5$ (and actually in arbitrary characteristic). It relies on Lemma\;\ref{5p5} below. Here we have to use the Kato--Milne cohomology group
 \[
 H^3(F):=H^3(F\,,\,\Z/2(2))
 \] for any field $F$ of characteristic 2. For basic facts about Kato--Milne cohomology, the readers are referred to \cite{Kato82} or \cite[Appendix]{GMS}.

\begin{lemma}[{See \cite[(5.1)]{KahnRostSujatha98} in characteristic $\neq 2$}]\label{5p5}
  When $\dim\vp\ge 5$ there is a natural isomorphism
  \[
  \theta\;:\;\ker\bigl(H^3(k)\lra H^3(k(X)\bigr)\simto \CH^2(X)_{\tors}\;.
  \]
\end{lemma}
\begin{proof}
To see that such an isomorphism exists one can simply apply \cite[Cor.\;7.1]{Kahn96DocMath}. For a more explicit construction, one can also
proceed as in the proof of \cite[\S\;2, Prop.\;1]{Merkurjev95KtheorySimpleAlg}.
\end{proof}

We need the following result, which is a characteristic 2 analogue of \cite[Satz\;5.6]{Arason75JA}.

\begin{thm}\label{5p6} Let $\varphi$ be a nondegenerate quadratic form of dimension $\ge 3$ over $k$ and $X=X_{\varphi}$ its projective quadric. Then, for every $\al\in H^3(k(X)/k):=\ker(H^3(k)\to H^3(k(X)))$, if $\al\neq 0$, there must exist elements $a,\,b\in k^*$ and $c\in k$ such that $\al=(a)\cup (b)\cup (c]\in H^3(k)$ and $\varphi$ is similar to a subform of the $3$-fold Pfister form $\Pfi{a,\,b\,;\,c}$.
\end{thm}
\begin{proof}
  Let $\psi$ be a 3-dimensional nondegenerate subform of $\varphi$. After scaling if necessary, we may assume $\psi=[1,\,c]\bot \langle b\rangle$. Write $k(\psi)$ for the function field of the projective quadratic defined by $\psi$ over $k$, and similarly for the function fields of other projective quadrics.

  Since $\varphi$ is isotropic over $k(\psi)$, the field extension $k(\psi)(\varphi)/k(\psi)$ is purely transcendental. Hence
  \[
  H^3\bigr(k(\psi)(\varphi)/k(\psi)\bigr):=\ker\left(H^3(k(\psi))\lra H^3\bigl(k(\psi)(\varphi)\bigr)\right)=0\;.
  \]Therefore, $H^3(k(\varphi)/k)\subseteq H^3(k(\psi)/k)$.

  From \cite[Thm.\;3.6]{AravireJacob09Contemp493} we know that $H^3(k(\psi)/k)=k^*\cup (b)\cup (c]$. In particular, every element $\al\in H^3(k(\varphi)/k)$ can be written as $\al=(a)\cup (b)\cup (c]$ for some $a\in k^*$. Let $\pi=\Pfi{a,\,b\,;\,c}$ be the Pfister form corresponding to $\al\in H^3(F)$. We assume $\al\neq 0$, so that $\pi$ is anisotropic over $k$.

The assumption $\al_{k(\varphi)}=0$ implies that $\pi_{k(\varphi)}\in I^4_q(k(\varphi))$ by \cite{Kato82Invent}. Using the Hauptsatz, we conclude that $\pi_{k(\varphi)}$ is hyperbolic. Then, it follows from \cite[Thm.\;4.2 (i)]{HoffmannLaghribi04TAMS} that $\varphi$ is similar to a subform of the 3-fold Pfister form $\pi$.
\end{proof}

\begin{remark}\label{5p7}
  In Theorem\;\ref{5p6} we can weaken the nondegeneracy assumption on $\vp$ to the following: \emph{$\vp$ is regular and not totally singular} (in the sense as defined in \cite[\S\;7, p.42]{EKM08}).

   Indeed, this weakened assumption still ensures that $\vp$ contains a 3-dimensional nondegenerate subform $\psi$, the field extension $k(\psi)(\varphi)/k(\psi)$ is still purely transcendental (cf. \cite[Prop.\;22.9]{EKM08}), and \cite[Thm.\;4.2 (i)]{HoffmannLaghribi04TAMS} is still valid.
\end{remark}

\begin{proof}[Proof of Theorem$\;\ref{5p3}$]
By Lemma\;\ref{5p5}, $\CH^2(X)_{\tors}$ is isomorphic to the kernel of the natural map
$\eta:H^3(k)\to H^3(k(X))$. If $\vp$ is isotropic, then $\ker(\eta)=0$. If $\vp$ is anisotropic, then by Thm.\;\ref{5p6}, $\ker(\eta)$ consists of symbols whose corresponding 3-fold Pfister form contains $\vp$ up to a scalar multiple. Since $\dim\vp\ge 5$, such a symbol is unique if it exists.
Thus, if $\eta$ is not injective, we have $\ker(\eta)\cong\Z/2$, and this case happens if and only if $\vp$ is an anisotropic neighbor of a 3-fold Pfister form. The theorem is thus proved.
\end{proof}

\section{Chow groups of affine quadrics}

To prepare the proofs of our results about codimension three Chow groups, we need some analysis on affine quadrics.

\medskip

We begin with a characteristic 2 variant of {\cite[(5.3)]{Karpenko95Nova}}.
\begin{lemma}\label{6p1}
  Let $\rho$ be a nondegenerate quadratic form of dimension $n\ge 2$ over $k$, and suppose that $\rho$ is not a hyperbolic plane. Let $a\in k$ and $\psi=\langle a\rangle\bot\rho$. Let $U\subseteq \A^n_k$ be the affine quadric defined by $a+\rho=0$.

  Then $\CH^p(U)=0$ in the each of the following cases:

  \begin{enumerate}
    \item The form $\psi$ is nondegenerate (i.e. $a\neq 0$ and $\dim\rho$ is even) and $\CH^p(X_{\psi})$ is elementary.
    \item $a=0$ and $\CH^p(X_{\rho})$ is elementary.
  \end{enumerate}
\end{lemma}
\begin{proof}
The proof in \cite[(5.3)]{Karpenko95Nova} works verbatim as soon as we notice that when $a=0$,  $\CH^p(X_{\psi})\cong \CH^p(X_{\rho})$ and the pushforward map $\CH^{p-1}(X_{\rho})\to \CH^p(X_{\psi})$ may be identified with the multiplication by $h\in \CH^1(X_{\rho})$ (\cite[Lemma\;70.2]{EKM08}).
\end{proof}

The following is easily deduced from  Lemma\;\ref{6p1},  Theorem\;\ref{5p3} and Prop.\;\ref{2p3}.
\begin{coro}[{Compare \cite[(5.4)]{Karpenko95Nova}}]\label{6p2}
   Let $\rho$ be an anisotropic (hence non-hyperbolic) nondegenerate quadratic form of dimension $n\ge 2$ over $k$. Let $a\in k$ and let $U\subseteq \A^n_k$ be the affine quadric defined by $a+\rho=0$.

   \begin{enumerate}
     \item Suppose $a=0$. Then $\CH^2(U)=0$ in the following cases:
   \begin{enumerate}

       \item $\dim\rho>8$.
       \item $\dim\rho\in\{5,\,7,\,8\}$, and $\rho$ is not a Pfister neighbor (e.g. $\rho$ is a $7$ or $8$ dimensional form  containing an Albert form).
       \item $\dim\rho=6$, and $\rho$ is neither an Albert form nor a Pfister neighbor.
   \end{enumerate}
     \item Suppose $a\neq 0$.  Then $\CH^2(U)=0$ in the following cases:

      \begin{enumerate}

       \item $\dim\rho$ is even and $\ge 8$.
       \item $\dim\rho=6$ and $\rho$ is not a Pfister neighbor.
       \item $\dim\rho=4$ and $\rho$ is not contained in a scalar multiple of a $3$-fold Pfister form

        \end{enumerate}
   \end{enumerate}
\end{coro}

\begin{lemma}\label{6p3}
Let $\rho$ be a nondegenerate quadratic form of dimension $n\ge 2$ over $k$, and suppose that $\rho$ is not a hyperbolic plane. Let $a\,,\,b\in k,\,c\in k^*$ and  $\vp=[ac^{-2}\,,\,b]\bot\rho$. Let $U\subseteq \A^{n+1}_k$ be the affine quadric defined by the equation $a+cy+by^2+\rho(x_1,\cdots, x_n)=0$.

If $\CH^p(X_{\vp})$ is elementary,  then $\CH^p(U)=0$.
\end{lemma}
\begin{proof}
  Let $\psi=\langle b\rangle \bot\rho$ (which can be degenerate). Note that $[ac^{-2},\,b]$ is isomorphic to the binary form $ax^2+cxy+by^2$. So we have the exact excision sequence
  \[
   \CH^{p-1}(X_{\psi})\overset{i_*}{\lra} \CH^p(X_{\vp})\lra \CH^p(U)\lra 0\,
  \]where the map $i_*$ is surjective when $\CH^p(X_{\vp})$ is elementary.
\end{proof}

\begin{coro}\label{6p4}
With notation and hypotheses as in Lemma$\;\ref{6p3}$, we have $\CH^2(U)=0$  in the following cases:

\begin{enumerate}
  \item  $\dim\rho>6$.
  \item $5\le \dim\rho\le 6$ and $\rho$ is not a Pfister neighbor.
\end{enumerate}
\end{coro}
\begin{proof}
 In the two cases above $\CH^2(X_{\vp})$ is elementary by Thm.\;\ref{5p3} and Prop.\;\ref{2p3}. Then apply Lemma\;\ref{6p3}.
\end{proof}

\section{Codimension three cycles on projective quadrics}\label{sec7}

In this section we prove our results about codimension three Chow groups.

\medskip

For a nondegenerate quadratic form $\vp$ over $k$, we write $\vp\in I^2_q(k)$ if $\dim\vp$ is even and $\mathrm{Arf}(\vp)=0$. If $\vp\in I_q^2(k)$ and $\vp$ has trivial Clifford invariant, we write $\vp\in I_q^3(k)$.

\begin{lemma}[{See \cite[XII.2.8]{Lam} in characteristic $\neq 2$}]\label{7p1}
  Let $\vp$ be a nondenegenerate quadratic form of dimension $10$ over $k$. If $\vp\in I^3_q(k)$, then $\vp$ is isotropic.
\end{lemma}
\begin{proof}
  We can write $\vp=\tau\bot\psi$ with $\psi$ nondegenerate of dimension 6. As $\vp$ has trivial Clifford invariant, the Brauer classes $[C(\psi)]$ and $[C(\tau)]$ coincide.

 If $\mathrm{Arf}(\psi)=0$, then $\mathrm{Arf}(\tau)=0$ and hence the 4-dimensional form $\tau$ is similar to a 2-fold Pfister form. It follows that the Brauer class $[C(\psi)]=[C(\tau)]$ has index $\le 2$. This implies that the Albert form $\psi$ is isotropic, and we are done.

Now we can assume $\mathrm{Arf}(\psi)\neq 0$ and $\psi$ is anisotropic. Let $K/k$ be the separable quadratic extension representing $\mathrm{Arf}(\psi)$. Then the above argument shows that $\psi_K$ is isotropic. By \cite[Prop.\;34.8]{EKM08}, there is a decomposition $\psi=a.N_{K/k}\bot\tau'$ for some $a\in k^*$ and some $4$-dimensional form $\tau'$. Since $\mathrm{Arf}(\psi)=\mathrm{Arf}(aN_{K/k})$, $\mathrm{Arf}(\tau')=0$  Setting $\psi'=\tau\bot a.N_{K/k}$, we are back to the situation $\vp=\tau'\bot\psi'$ with $\psi'$ an Albert form. The argument in the previous paragraph shows that $\psi'$ is isotropic. The lemma is thus proved.
\end{proof}

\begin{lemma}\label{7p2}
Let $\vp$ be a nondegenerate quadratic form of dimension $10$ over $k$. Suppose that $\vp\in I^2_q(k)\setminus I^3_q(k)$.

Then there exists an odd degree extension $K/k$ and a separable extension $L/K$ with $[L:K]=2^{5-s}$ such that  $\vp_L$ is hyperbolic, where $s=s(\vp)$ is the splitting index of $\vp$.
\end{lemma}
\begin{proof}
  By \eqref{3p1p1}, the assumption $\vp\notin I^3_q(k)$ means that  $C(\vp)$ does not split, whence $s(\vp)\le 3$.

  First assume $s(\vp)=3$. Let $F/k$ be a separable quadratic extension such that some binary nondegenerate subform of $\vp$ becomes isotropic over $F$. Then $\vp_F=\bH\bot\rho_F$ for some 8-dimensional form $\rho_F\in I^2_q(F)$. Then $s(\rho_F)=s(\vp_F)-1\ge s(\vp)-1=2$. By Lemma\;\ref{3p5} (1), we can find a quadratic separable extension $L/F$ such that $\rho_L$ is hyperbolic. Now $[L:k]=4=2^{5-s}$ and we can take $K=k$.

  Now let us assume $s=s(\vp)\le 2$. By \cite[\S\,15.2, Lemma]{PierceGTM88}, there exists an odd degree extension $K/k$ and a separable extension $F/K$ of degree $2^{3-s}$ such that $\ind(C(\vp)_F)=2$. Then $s(\vp_F)=3$. So by the previous case we can find a separable  extension $L/F$ of degree 4 such that $\vp_L$ is hyperbolic. Now $[L:K]=2^{3-s}\cdot 4=2^{5-s}$. The lemma is thus proved.
\end{proof}

\begin{thm}\label{7p3}
  Let $X=X_{\vp}$ be the projective quadric defined by a nondegenerate quadric form $\vp$ over $k$.

Then $\big|\CH^3(X)_{\tors}\big|\le 2$.
\end{thm}
\begin{proof}
  If $\vp$ is isotropic, then $\CH^3(X)_{\tors}\cong \CH^{2}(Y)_{\tors}$ for a lower dimensional smooth quadric $Y$. In this case the theorem follows from the results for Chow groups of codimension $2$ (Theorem\;\ref{5p3}).

  Now we can assume $\vp$ is anisotropic. Note that $\CH^3(X)\cong K_0(X)^{(3/4)}$ (cf. \eqref{4p1}). If $\vp\notin I^2_q(k)$, we can just apply \cite[(3.8)]{Karpenko90AlgGeoInv}. So we assume $\vp\in I^2_q(k)$. In particular $\dim\vp$ is even.

  If $\dim\vp\le 8$, i.e., $m:=\frac{\dim X}{2}\le 3$, then $2m-3\le m$. With notation as in \cite[(3.10)]{Karpenko90AlgGeoInv}, in the torsion subgroup of the second kind the dimension $2m-3$ component $T^{II}_{2m-3}$ is 0 and hence
  \[
  \CH^3(X)_{\tors}\cong \left(K_0(X)^{(3/4)}\right)_{\tors}=T^I_{2m-3}\cong \Z/2\text{ or } 0\,.
  \]It remains to consider the case where $\vp\in I^2_q(k)$, $\dim\vp\ge 10$ and $\vp$ is anisotropic.

  Now $K_0(X)^{(i/i+1)}\cong \CH^i(X)$ is torsion free for $i\le 2$. (For $i=2$ we use Thm.\;\ref{5p3}.) By the last assertion in \cite[(3.9)]{Karpenko96DocMath}, $(T^I)^{(3)}=0$ and hence
  $\CH^3(X)_{\tors}\cong \left(K_0(X)^{(3/4)}\right)_{\tors}=(T^{II})^{(3)}$ is a cyclic group.

  It is now sufficient to show that $\CH^3(X)_{\tors}$ is killed by 2.

  If $\dim\vp>10$, then we can write $\vp=\rho\bot\tau$ with $\dim\tau=2$ and  $\dim\rho>8$. Choosing $L/k$ to be a quadratic separable extension with $\tau_L\cong\bH$, we get
  $\CH^3(X_L)_{\tors}\cong \CH^2(Y_L)_{\tors}$, where $Y$ is the quadric defined by $\rho$.
 Here $\CH^2(Y_L)_{\tors}=0$ by Thm.\;\ref{5p3}. So the standard restriction-corestriction argument shows that $2\cdot\CH^3(X)_{\tors}=0$.

So now we assume $\dim\vp=10$ (and $\vp$ is anisotropic, belonging to $I^2_q(k)$).

Since $\vp$ is anisotropic, Lemma\;\ref{7p1} implies  $\vp\notin I^3_q(k)$. Let $s=s(\vp)$. By Lemma\;\ref{7p2}, we can find an odd degree extension $K/k$ and a separable extension $L/K$ of degree $2^{5-s}$ such that $\vp_L$ is hyperbolic. Note that the splitting index does not change after  an odd degree base extension. So $s(\vp_K)=s(\vp)=s$. Now, by the estimate of $|T^{II}|$ in \cite[(3.9)]{Karpenko96DocMath} we have
\[
\big|\CH^3(X)_{\tors}\big|\le \big|\CH^3(X_K)_{\tors}\big|=\Big|\left(T^{II}\right)^{(3)}\Big|\le |T^{II}|\le 2^{s+(5-s)-4}=2\,.
\]The theorem is thus proved.
\end{proof}

\begin{remark}\label{7p4}
Our proof of Theorem\;\ref{7p3} is slightly different from Karpenko's arguments (\cite[\S\,3]{Karpenko91cyclesCodim3} or \cite[\S\,5]{Karpenko96DocMath}). His approach relied on the following result (which is part of the second assertion in Theorem\;\ref{5p3}):

\emph{For a 8-dimensional form $\rho$, $\CH^2(X_{\rho})_{\tors}\neq 0$ if and only if $\rho$ is similar to a 3-fold Pfister form.}

Our proof here does not need any characterization of 8-dimensional forms with nontrivial torsion in the codimension 2 Chow group.
 We have only used the first assertion of Theorem\;\ref{5p3} and the vanishing of $\CH^2(X_{\vp})_{\tors}$ for $\vp$ of dimension $>10$. These two results can be proved without using Lemma\;\ref{5p5} (cf. Remark\;\ref{5p4}).
\end{remark}

We now prove  \cite[(6.2)]{Karpenko95Nova} in characteristic 2.

\begin{lemma}\label{7p5}
  Let $p,\,n\in\N$ with $n>2p+2$. Let $\mathcal{P}(p,\,n)$ be the following statement: For every extension field $F$ of $k$ and every nondegenerate quadratic form $\psi$ of dimension $n$ over $F$, the group $\CH^p(X_{\psi})$ is elementary.

  Then $\mathcal{P}(p,\,n)$ implies $\mathcal{P}(p,\,n+1)$.
\end{lemma}
\begin{proof}
  It is clear that $\mathcal{P}(0,\,n)$ holds for all $n>2$. We may thus assume $p\ge 1$.

  Let $F$ be a field extension of $k$ and let $\rho$ be a nondegenerate quadratic form of dimension $n-2$ over $F$. Then $\psi=\rho\bot\bH$ has dimension $n$ and $\CH^p(X_{\psi})\cong \CH^{p-1}(X_{\rho})$. So $\mathcal{P}(p,\,n)$ implies $\mathcal{P}(p-1,\,n-2)$, and by induction on $p$ we find that $\mathcal{P}(p,\,n)$ implies $\mathcal{P}(p-1,\,N)$ for all $N\ge n-2$.

  \medskip

  Now suppose $\mathcal{P}(p,\,n)$ holds and consider a nondegenerate quadratic form $\vp$ of dimension $n+1$ over $F$.   We distinguish two cases to show $\CH^p(X_{\vp})$ is elementary.

  \medskip

  \noindent{\bf Case 1.} $n+1$ is even.

  In this case we can write $\vp=[a,\,b]\bot\rho$ for some $(n-1)$-dimensional form $\rho$ over $F$. Let $U$ be the affine quadric defined by $a+y+by^2+\rho=0$. By induction and excision, it remains to show $\CH^p(U)=0$. This can be done by using a fibration arugment as in \cite[(6.2)]{Karpenko95Nova} and applying Lemma\;\ref{6p1}.

  \medskip

  \noindent{\bf Case 2.} $n+1$ is odd.

  Now we can write $\vp=\langle a\rangle\bot [b,\,c]\bot\tau$ for some nondegenerate form $\tau$ of dimension $n-2$.  If $b=0$, then  $[b,\,c]\cong \bH$ and $\CH^p(X_{\vp})\cong\CH^{p-1}(X_{\langle a\rangle\bot\tau})$.  The result then follows immediately from $\mathcal{P}(p-1,\,n-1)$.

  We may thus assume $b\neq 0$.   Let $U$ be the affine quadric defined by
  \[
  a+bx^2+xy+cy^2+\tau=0\,.
  \]
  As in Case 1, it is sufficient to show $\CH^p(U)=0$. This again follows from the induction hypothesis and a fibration argument, together with Lemma\;\ref{6p1}.
\end{proof}

\begin{prop}\label{7p6}
  Let $n$ be an odd integer $>8$. Then the following are equivalent:

  \begin{enumerate}
    \item For every field extension $F/k$ and every nondengenerate quadratic form $\psi$ of dimension $\ge n$ over $F$, $\CH^3(X_{\psi})_{\tors}=0$.
    \item For every field extension $F/k$ and every nondengenerate quadratic form $\psi$ of dimension $n$ over $F$, $\CH^3(X_{\psi})_{\tors}=0$.
    \item For every field extension $F/k$ and every nondengenerate quadratic form $\psi$ of dimension $n+1$ over $F$ with $\psi\in I^2_q(F)$, $\CH^3(X_{\psi})_{\tors}=0$.
    \item For every field extension $F/k$ and every nondengenerate quadratic form $\psi$ of dimension $n+1$ over $F$ with $\psi\in I^3_q(F)$, $\CH^3(X_{\psi})_{\tors}=0$.
  \end{enumerate}
\end{prop}
\begin{proof}
  Combine Lemma\;\ref{7p5}, Prop.\;\ref{4p2} and \cite[(4.9)]{Karpenko95Nova}. (The proof of the result we cite from \cite{Karpenko95Nova} goes through in characteristic 2 without change.)
\end{proof}

\begin{lemma}\label{7p7}
  Let $\tau$ be a nondegenerate quadratic form of even dimension $m\ge 6$ over $k$, and  let $U_0\subseteq\A^{m+3}$ be the affine quadric over $k$ defined by the equation
  \[
  a_1+cY_1+b_1Y_1^2+[a_2,\,b_2]\bot\tau=0\;,\quad \text{where }\; c\,,\,a_i,\,b_i\in k^* \,.
  \]

  Assume either $m\ge 8$ or $\tau$ is an Albert form.

  Then  $\CH^3(U_0)\cong \CH^3(U)$, where $U\subseteq\A^{m+1}$ is the affine quadric over the rational function field $F=k(y_1,\,x_2)$ defined by the equation
  \[
  (a_1+cy_1+b_1y_1^2+a_2x_2^2)+x_2Y_2+b_2Y_2^2+\tau=0\,.
  \]
\end{lemma}
\begin{proof}
   Let $k_1$ be the rational function field $k(y_1)$ and let $U_1\subseteq\A^{m+2}$ be the affine quadric over $k_1$ defined by
  \[
  (a_1+cy_1+b_1y^2_1)+[a_2,\,b_2]\bot\tau=0\,.
  \]
  By considering a fibration over $\A^1$ as in Case 1 of the proof of Lemma\;\ref{7p5}, we can use Cor.\;\ref{6p2} to get $\CH^3(U_0)\cong \CH^3(U_1)$.

  Now letting $x_2,\,y_2$ denote the variables corresponding to the binary form $[a_2,\,b_2]$, consider the projection $\pi: U_1\to \A^1_{k_1}$ onto the $x_2$-coordinate. Then the generic fiber of $\pi$ is the affine quadric $U$ over $F=k_1(x_2)=k(y_1,\,x_2)$  in the statement of the lemma.
 By the fibration method, to show  $\CH^3(U_1)\cong\CH^3(U)$, it is sufficient to prove that for every closed point $P\in \A^1_{k_1}$, the closed fiber $(U_1)_P$ of $\pi$ over $P$ satisfies $\CH^2((U_1)_P)=0$.

 Let us fix a closed point $P\in \A^1_{k_1}$ and put $V=(U_1)_P$.  Writing $\al=a_1+cy_1+b_1y_1^2\in k_1$, $V$ is the affine quadric over $K:=k_1(P)$ defined by the equation
 \[
 (\al+a_2x_2(P)^2)+x_2(P)Y_2+b_2Y_2^2+\tau=0\,.
 \]
 If $x_2(P)\neq 0$, we can deduce from Cor.\;\ref{6p4} that $\CH^2(V)=0$. If $x_2(P)=0$, then $V$ is defined by $\al+\langle b_2\rangle\bot\tau=0$. In this case, we have a fibration exact sequence
 \[
 \bigoplus_{Q\in \A^1_K}\CH^1(V_Q)\lra \CH^2(V)\lra \CH^2(V_{\eta})\lra 0\,,
 \]where the generic fiber $V_{\eta}$ is the affine quadric defined by $(\al+b_2y_2^2)+\tau=0$ over the rational function field $K(y_2)$. By Cor.\;\ref{6p2} (2), we have $\CH^2(V_{\eta})=0$. For each closed point $Q\in\A^1_K$, Lemma\;\ref{6p1} shows that $\CH^1(V_Q)=0$. So we get $\CH^2(V)=0$ as desired. The lemma is thus proved.
\end{proof}

\begin{lemma}\label{7p8}
  Let $\vp$ be a nondegenerate $14$-dimensional quadratic form over $k$. Suppose that $\vp$ contains an Albert form as a subform.

  Then $\CH^3(X_{\vp})$ is elementary.
\end{lemma}
\begin{proof}
  We may assume $\vp$ is anisotropic and write $\vp=[a_1,\,b_1]\bot\cdots\bot[a_4,\,b_4]\bot\rho$, where $a_i,\,b_i\in k^*$ and $\rho$ is an Albert form.
  Put $F=k(y_1,\,y_2,\,x_2,\,x_3)$ and let $U$ be the affine quadric over $F$  defined by the equation
  \[
  (\al+a_4x_4^2)+x_4Y_4+b_4Y^2_4+\rho=0\;\; \text{ where } \al=a_1+y_1+b_1y_1^2+\sum_{2\le i\le 3}(a_ix_i^2+x_iy_i+b_iy_i^2)\;.
  \] By using the standard excision sequence, a repeated application of Lemma\;\ref{7p7} shows $\CH^3(X_{\vp})/\Z.h^3\cong \CH^3(U)$.

By Lemma\;\ref{6p3}, it suffices to show $\CH^3(X_{\theta})$ is elementary, where $\theta$ is the form $[\al x_4^{-2}+a_4,\,b_4]\bot\rho$ over $F$.

  Put $\psi=\langle b_4\rangle\bot\rho$. Then $\psi$ is anisotropic since $\vp$ is, and it is not a Pfister neighbor since the Albert form $\rho$ is not a Pfister neighbor.
  Therefore, $\CH^2(X_{\psi})$ is elementary by Thm.\;\ref{5p3}. Now the groups $K_0(X_{\psi})^{(i/i+1)}\cong \CH^i(X_{\psi}),\,i\le 2$ are elementary. Therefore, using Prop.\;\ref{4p3} we reduce the problem to proving the following assertion: the discriminant algebra $K$ of the form $\theta$ over $F$ is a field such that $\ind(C_0(\psi)_K)=\ind(C_0(\psi))$.

      In fact, $K$ is the function field $k(\tau)$ of the quadratic form $\tau$ over $k$ given by
      \[
      \tau=b_4([a_1,\,b_1]\bot[a_2,\,b_2]\bot [a_3,\,b_3])\bot [1,\,a_4b_4]\;.
      \]Note that $C_0(\psi)$ is not a division algebra since the Albert form $\rho$ is a subform of $\psi$. So the division algebra Brauer equivalent to $C_0(\psi)$ has dimension $<\dim C_0(\psi)=2^6$. Since $\dim\tau=8\ge 6+1$, it follows from Merkurjev's index reduction theorem (\cite[Cor.\;30.9]{EKM08}) that $\ind(C_0(\psi)_{K})=\ind(C_0(\psi)_{k(\tau)})=\ind(C_0(\psi))$. This completes the proof of the lemma.
\end{proof}

\newpara\label{7p9} We recall some facts on residue forms in the case of valued fields. Let $A$ be a ring endowed with a rank $1$ discrete Henselian valuation $\nu$. Let $K$ and $A^{\times}$ be its field of fractions and the group of units, respectively. Let $\pi$ be a uniformizing parameter and $k=A/\pi A$ the residue field. Let $\vp$ be an anisotropic quadratic form over a $K$-vector space $V$. Since $\vp$ is anisotropic and $\nu$ is Henselian, we have the following inequality:
\begin{equation}\label{7p9p1}
\nu(B_{\vp}(x,y)^2)\geq \nu(\vp(x))+\nu(\vp(y))
\end{equation}
for all $x, y\in V$ (\cite[Lemma 2.2]{Tietze}).

For $i \in \Z$, let $V_i=\{x\in V\mid \vp(x)\in \pi^i A\}$. Using the inequality (\ref{7p9p1}), we prove that $V_i$ is an $A$-module. The form $\vp$ induces two quadratic forms $\overline{\vp_0}$ and $\overline{\vp_1}$, called the first and the second residue forms, on the $k$-vector space $V_i/V_{i+1}$ as follows:
\begin{eqnarray*}
\overline{\vp_i}: V_i/V_{i+1} & \longrightarrow & k\\ x+V_{i+1}& \mapsto & \overline{\pi^{-i}\vp(x)}
\end{eqnarray*}

Obviously, the quadratic forms $\overline{\vp_0}$ and $\overline{\vp_1}$ are anisotropic. When $\vp$ is nonsingular, we have by \cite[Theorem 1]{Mammone-Moresi-Wadsworth}:
\begin{equation}\label{7p9p2}
\dim \vp=\dim \overline{\vp_0} +\dim \overline{\vp_1}.
\end{equation}

The following result  is an analogue in characteristic $2$ of a theorem of Rost (cf. \cite{Rost99dim14}, \cite{Rost06}).

\begin{thm}\label{7p10}
Let $\vp \in I^3_q(k)$ be an anisotropic form of dimension $14$. Then, $\vp$ contains an Albert form as a subform.
\end{thm}
\begin{proof}
Let $A$ be a Henselian discrete rank $1$ valuation ring of characteristic $0$ whose maximal ideal is $2A$ and residue field $k$ (see \cite[(1.4)]{Wadsworth}). Let $K$ and $A^{\times}$ be the field of fractions and the group of units of $A$, respectively.

There exists a nondegenerate quadratic module $\theta$ of rank $14$ defined on an $A$-module $V$ that is a lifting of $\vp$, i.e., $\vp$ is isometric to the induced quadratic form $\overline{\theta}$ on the $k$-vector space $V/2V$. The form $\theta$ is anisotropic.

Let $S= \{(-1)^ka^2+4b \mid k\in \Z, a\in A^{\times}, b\in A\}$. This is clearly a subgroup of $A^{\times}$. By \cite[Lemma 1.6]{Wadsworth}, there exists a surjective group homomorphism $$\gamma: S \longrightarrow k/\wp(k)$$ given by: $(-1)^ka^2+4b\mapsto ba^{-2} +\wp(k)$, and ${\rm Ker}\gamma=\pm A^{\times 2}$. Moreover, $\det \theta \in S/A^{\times 2}$ and $\gamma(\det \theta)= \mathrm{Arf}(\overline{\theta})$ \cite[Proposition 1.14]{Wadsworth}.

Using \cite[Corollary 5.4]{BarryChapmanLaghribi20Israel}, we get a form $\vp'\in I^3A$ such that $\overline{\vp'}$ is Witt-equivalent to $\vp \cong \overline{\theta}$. Since $A$ is Henselian, it follows that $\vp'$ is Witt-equivalent to $\theta$ \cite[Satz 3.3]{Knebusch69}. Hence, $\theta \in I^3A$. In particular, $\theta_K\in I^3K$. It follows from a theorem of Rost (\cite{Rost99dim14}, \cite{Rost06}) that $\theta_K$ contains an Albert form $\theta'$ as a subform.

We write $\theta'\cong [a_1, b_1] \perp [a_2, b_2] \perp [a_3, b_3]$ for suitable $a_i, b_i \in K$, $1\leq i\leq 3$. We claim that $a_i, b_i \in A^{\times}$ for all $1\leq i\leq 3$, i.e., $\theta'$ is defined over $A$. In fact, let us write $a_i=u_i2^{\epsilon_i}$ and $b_i=v_i2^{\delta_i}$ for $u_i, v_i \in A^{\times}$ and
$\epsilon_i, \delta_i\in \Z$.

(i) The form $\vp$ is nothing but the first residue form of $\theta_K$, and thus the second residue form of $\theta_K$ is the zero form by (\ref{7p9p2}).

(ii) By (i) we deduce that $\epsilon_i$ and $\delta_i$ are even for all $1\leq i\leq 3$, otherwise the second residue form of $\theta_K$ would be of dimension $>0$.

(iii) By (ii) we get $[a_i, b_i]\cong [u_i2^{\epsilon_i + \delta_i}, v_i]$ for all $1\leq i\leq 3$ (using that $\delta_i$ is even and the isometry $a[b, c]\cong [ab, a^{-1}c]$ for scalars $a\neq 0, b$ and $c$).

(iv) By the inequality (\ref{7p9p1}), we have $\epsilon_i + \delta_i\leq 0$. Moreover, if for some $i$, we have $\epsilon_i + \delta_i<0$, then the first residue form of $[u_i2^{\epsilon_i + \delta_i}, v_i]$ is the degenerate form $\left<\overline{u_i}, \overline{v_i}\right>$, this is excluded since $\vp$ is nondegenerate. Consequently, $[a_i, b_i]=[u_i2^{\epsilon_i}, v_i2^{-\epsilon_i}]\cong 2^{\epsilon_i}[u_i, v_i] \cong [u_i, v_i]$ because $\epsilon_i$ is even.

Hence, $\theta'\cong (\theta'')_K$, where $\theta''$ is the form $[u_1, v_1]\perp [u_2, v_2] \perp [u_3, v_3]$.

Now, the conditions that $\theta''$ is defined over $A$ and $\theta\cong (\theta'')_K$ is a subform of $\theta_K$ imply that $\theta''$ is also a subform of $\theta$ over $A$. Taking the reduction modulo $2$, we get that $\overline{\theta''}$ is a subform of $\vp$. The form $\theta''$ has determinant $-A^{\times 2}$ because the scalar $\prod\limits_{1\leq i\leq 3}(4u_iv_i-1)\in A^{\times}$ is a representative of $\det (\theta'')_K= -K^{*2}\in K^*/K^{*2}$. Since ${\rm Ker}\gamma= \pm A^{\times 2}$, it follows that $\gamma(\det \theta'')= \mathrm{Arf}(\overline{\theta''})=0$, which means that $\overline{\theta''}$ is an Albert form.
\end{proof}

\begin{thm}\label{7p11}
  For every nondegenerate form $\vp$ of dimension $\ge 13$ over $k$, $\CH^3(X_{\vp})$ is elementary.
\end{thm}
\begin{proof}
  Combine Theorem\;\ref{7p10} with Lemma\;\ref{7p8} and Prop.\;\ref{7p6}.
\end{proof}

In characteristic different from 2, Izhboldin completely determined when the group $\CH^3(X_{\vp})_{\tors}$ is trivial for all nondegenerate forms $\vp$ of dimension $\ge 9$ (\cite[Thm.\;0.5]{Izhboldin01}). A full proof of his theorem builds upon computations of the fourth unramified cohomology groups of quadrics. Without going into study of unramified cohomology, in the rest of this section we discuss a few cases of Izhboldin's results in characteristic 2.

\

We begin with the following result, which is the characteristic 2 analogue of \cite[Prop.\;3.7]{Izhboldin01}.

\begin{prop}\label{7p12}
Let $\varphi$ be a nondegenerate quadratic form over $k$ satisfying one of the following conditions:

\begin{enumerate}
  \item $\dim\varphi=12$, $\mathrm{Arf}(\varphi)\neq 0\in k/\wp(k)$, and $\ind(\varphi)\le 2$.
  \item $\dim\varphi=11$ and $\ind(\varphi)\ge 2$.
  \item $\dim\varphi=10$, $\mathrm{Arf}(\varphi)\neq 0\in k/\wp(k)$, and $\ind(\varphi)=2$.
  \item $\dim\varphi=9$ and $\ind(\varphi)\ge 4$.
\end{enumerate}

Then $\CH^3(X_{\varphi})_{\tors}=0$.
\end{prop}

Let us consider the subcase with $\ind(\varphi)=2$ in Prop.\;\ref{7p12} (1). By Theorem\;\ref{7p11}, we have $\CH^3(X_{\tau})_{\tors}=0$ for every nondegenerate form $\tau$ of dimension $15$. Applying  \cite[Prop.\;2.8 (iii)]{Izhboldin01} with $n=7$ (which is still valid in characteristic 2, since its proof works almost verbatim), we see that $\CH^3(X_{\varphi})_{\tors}=0$.

The proofs of the other cases in Prop.\;\ref{7p12}  may also be given along the lines of the case treated in \cite[Prop.\;3.7]{Izhboldin01}. We shall not provide full details, but content ourselves with the easy observation that the key ingredient we need is the characteristic 2 version of \cite[Lemma\;1.19]{Izhboldin01}. That is, it suffices to prove the following:

\begin{lemma}\label{7p13}
  Let $n$ be an integer $\ge 2$ and let $\varphi$ be a nondegenerate quadratic form over $k$ such that one of the following conditions holds:

\begin{enumerate}
  \item $\dim\varphi=2n$, $\mathrm{Arf}(\varphi)\neq 0\in k/\wp(k)$, and $\ind(\varphi)\le 2$.
  \item $\dim\varphi=2n-1$ and $\ind(\varphi)\ge 2$.
  \item $\dim\varphi=2n-2$, $\mathrm{Arf}(\varphi)\neq 0\in k/\wp(k)$, and $\ind(\varphi)=2$.
  \item $\dim\varphi=2n-3$ and $\ind(\varphi)\ge 4$.
\end{enumerate}

Then there exists a $(2n+1)$-dimensional nondegenerate form $\tilde{\varphi}$ and a $(2n+2)$-dimensional nondegenerate form $\gamma\in I^3_q(k)$ such that $\varphi\subseteq\tilde\varphi\subseteq\gamma$ and $\ind(\tilde\varphi)=1$.
\end{lemma}

Below we provide a detailed the proof of Lemma\;\ref{7p13}, which seems to involve some more subtleties than its counterpart in characteristic different from 2.

\medskip

First note that we have:

\begin{lemma}\label{7p14}
  Let $A$ be a central simple $k$-algebra of exponent $\le 2$, $L/k$ a separable field extension of degree $\le 2$ and $m$ an integer. Suppose that one of the following conditions holds:

  \begin{enumerate}
    \item $\ind(A_L)=1$ and $m=2$.
    \item $L=k$, $\ind(A)\le 2$, and $m=3$.
    \item $\ind(A_L)\le 2$ and $m=4$.
    \item $L=k$, $\ind(A)\le 4$, and $m=5$.
  \end{enumerate}

  Then there exists an $m$-dimensional nondegenerate form $\mu$ over $k$ such that the algebra $C'_0(\mu)$ has center $L$ and is Brauer equivalent to $A_L$.
\end{lemma}
\begin{proof}
  In Cases (1)--(3), one can use the same arguments as in the proof of \cite[Lemma\;1.17]{Izhboldin01}. It suffices to change the notations
  \[
  k(\sqrt{d})\,,\;\pfi{d}\,,\; \dgf{1,\,-a,\,-b}\,,\;\pfi{a,\,b} \quad \text{in characteristic }\neq 2
  \]to
  \[
  k[t]/(t^2-t-d)\,,\;\Pfi{d}\,,\;\dgf{a}\bot [1,\,b]\,,\; \Pfi{a\,;\,b}\quad \text{ in characteristic }2\,.
  \]
  In Case (4), $A$ is Brauer equivalent to a biquaternion $k$-algebra, which gives rise to an Albert form $q=c[1,\,a]\bot\rho$, where $c\in k^*$, $a\in k$ and $\rho$ is a 4-dimensional form with Arf invariant $a\in k/\wp(k)$. Then we can take $\mu=\dgf{c}\bot \rho$.
\end{proof}

The proof of Lemma\;\ref{7p13} also relies on the lemma below.

\begin{lemma}[{See \cite[Lemma\;4.3]{Izhboldin98isotropy} in characteristic $\neq 2$}]\label{7p15}
  Let $\vp$ and $\psi$ be even dimensional nondegenerate quadratic forms over $k$ with the same Arf invariant (so the algebras $C'_0(\vp)$ and $C'_0(\psi)$ have the same center). Suppose that $C'_0(\vp)$ and $C'_0(\psi)$ are Brauer equivalent.

  Then there exists an element $a\in k^*$ such that $\vp\bot a\psi\in I^3_q(k)$.
\end{lemma}
\begin{proof}
  First assume $\vp$ and $\psi$ have trivial Arf invariant, i.e. they lie in $I^2_q(k)$. Then the assumption implies that $\vp$ and $\psi$ have the same Clifford invariant. So we can just take $a=-1$.

  Now assume $\vp$ and $\psi$ have the same nontrivial Arf invariant $d\in k/\wp(k)$. Then their discriminant algebra $L=k[t]/(t^2-t-d)$ is a quadratic separable field extension of $k$. By assumption in the Brauer group $\Br(L)$ we have
  \[
  [C(\vp)_L]=[C_0(\vp)]=[C_0(\psi)]=[C(\psi)_L]\,.
  \]Thus the forms $\tilde{\vp}:=\vp\bot [1,\,d]$ and $\tilde{\psi}=\psi\bot [1,\,d]$ lie in $I^2_q(k)$, and the forms $\tilde{\vp}_L$ and $\tilde{\psi}_L$ have the same Clifford invariant (see e.g. \cite[\S\;14]{EKM08}). Thus, the Clifford invariant $e_2(\tilde{\vp}-\tilde{\psi})$ of $\tilde{\vp}-\tilde{\psi}$ lies in
  \[
  \Br(L/k)=\ker(\Br(k)\lra \Br(L))\,.
  \]By the well known structure of the group $\Br(L/k)$, we have $e_2(\tilde{\vp})-e_2(\tilde{\psi})=(a,\,d]$ for some $a\in k^*$. Note that $e_2(a\tilde{\psi})=e_2(\tilde{\psi})$ since $\tilde{\psi}\in I^2_q(k)$. Thus $e_2(\tilde{\vp}-a\tilde{\psi}-\Pfi{a\,;\,d})=0$, and it follows that
  \[
  \vp-a\psi=\vp-a\psi+[1,\,d]-a[1,\,d]-\Pfi{a\,;\,d}=\tilde{\vp}-a\tilde{\psi}-\Pfi{a\,;\,d}\;\in \;I^3_q(k)\;.
  \]This completes the proof.
\end{proof}

\begin{proof}[Proof of Lemma$\;\ref{7p13}$]
Put $m=2n+2-\dim\varphi$.

In Cases (1) and (3), if $\varphi\in I^2_q(k)$, we put $A=C'_0(\varphi)$ and $L=k$; otherwise  put $A=C(\varphi)$ and let $L$ be the center of $C'_0(\varphi)=C_0(\varphi)$.  Then $A_L$ is Brauer equivalent to $C'_0(\varphi)$, and $\ind(\varphi)=\ind(A_L)$.  By Cases (1) and (3) of Lemma\;\ref{7p14}, there exists an $m$-dimensional nondegenerate form $\mu$ over $k$  such that $C'_0(\varphi)$ and $C'_0(\mu)$ have the same center and are Brauer equivalent. Here $\varphi$ and $\mu$ are even dimensional. So we can apply Lemma\;\ref{7p15} to find an element $a\in k^*$ such that the form $\gamma:=\varphi\bot a\mu$ lies in $I^3_q(k)$. Writing $a\mu=\theta\bot c[1,\,b]$ and setting $\tilde{\varphi}=\varphi\bot \theta\bot \dgf{c}$, we have
\[
[C_0(\tilde{\varphi})]=[C(c\gamma)]=[C(\gamma)]=0 \text{ in }\Br(k)
\]whence $\ind(\tilde{\varphi})=1$. Thus the forms $\gamma$ and $\tilde{\varphi}$ have the required properties, and we obtain the desired result.

Now consider Cases (2) and (4). We put $A=C'_0(\varphi)$ and $L=k$. By  Cases (2) and (4) of Lemma\;\ref{7p14}, there exists an $m$-dimensional nondegenerate form $\mu$ over $k$  such that $C'_0(\varphi)=C_0(\varphi)$ and $C'_0(\mu)=C_0(\mu)$ are Brauer equivalent over $k$. Write $\varphi=\rho_0\bot \dgf{a}$, $\mu=\mu_0\bot \dgf{ab}$ and choose representatives $a_0,\,b_0\in k$ of the Arf invariants $\mathrm{Arf}(\varphi_0),\,\mathrm{Arf}(\mu_0)$. Set
\[
\varphi_1:=\varphi_0\bot a[1,\,a_0+u_0]\,.
\]Then $[C(a\varphi_1)]=[C_0(\varphi)]$ and $[C(ab\mu_0)]=[C_0(\mu)]$ in $\Br(k)$. Also, it is easy to see that the form $\gamma:=\varphi_1\bot b\mu_0$ has trivial Arf invariant, i.e., $\gamma\in I^2_q(k)$. Now
\[
e_2(a\gamma)=e_2(a\varphi_1)-e_2(ab\mu_0)=[C_0(\varphi)]-[C_0(\mu)]=0\;\in \;\Br(k)\;.
\]It follows that $a\gamma\in I^3_q(k)$ and hence $\gamma\in I^3_q(k)$.

Set $\tilde{\varphi}:=\varphi_0\bot \dgf{a}\bot b\mu_0=\varphi\bot b\mu_0$. We have
\[
\varphi\subseteq\tilde{\varphi}=\varphi\bot b\mu_0=\varphi_0\bot \dgf{a}\bot b\mu_0\subseteq \gamma=\varphi_0\bot a[1,\,a_0+b_0]\bot b\mu_0\,,
\]and
\[
0=[C(a\gamma)]=[C_0(\varphi_0\bot \dgf{a}\bot b\mu_0)]=[C_0(\tilde{\varphi})]\;\in\;\Br(k)\,.
\]Hence $\ind(\tilde{\varphi})=1$. This completes the proof.
\end{proof}

\begin{remark}\label{7p16}
One can also check that Corollary\;3.10 and Lemmas\;3.11 and 3.12 of \cite{Izhboldin01} extend to characteristic 2. Namely, for a nondegenerate quadratic form $\varphi$ over $k$, the following statements hold:

\begin{enumerate}
  \item  Suppose $\dim\varphi$ is even and $>8$, the discriminant algebra $L$ of $\varphi$ is a field (i.e. $\mathrm{Arf}(\varphi)\neq 0$) and $\varphi_L$ is hyperbolic. Then $\CH^3(X_{\varphi})_{\tors}=0$.
  \item Suppose $\dim\varphi=10$, the discriminant algebra $L$ of $\varphi$ is a field (i.e. $\mathrm{Arf}(\varphi)\neq 0$) and $\varphi=\tau\bot c.N_{L/k}$ for some $c\in k^*$ and some subform $\tau$. Then  $\CH^3(X_{\varphi})_{\tors}=0$ except possibly when the following conditions hold simultaneously:
      \[
      \ind(\varphi)=\ind(\tau_L)=1\,,\;\ind(\tau)=2\quad\text{and}\quad \varphi_L \text{ is not hyperbolic}.
      \]
  \item Suppose $\dim\varphi=9$, $\ind(\varphi)>1$ and $\varphi$ has one of the following forms:

  (i) $\varphi=\gamma\bot [a,\,b]$, where $a,\,b\in k$ and $\gamma$ is a 7-dimensional Pfister neighbor.

  (ii) $\varphi=\tau\bot \dgf{d}$, where $d\in k^*$ and $\tau\in I^2_q(k)$.

  Then $\CH^3(X_{\varphi})_{\tors}=0$.
\end{enumerate}

Indeed, the above assertions follow on parallel lines along the proofs of the corresponding results in \cite{Izhboldin01}, as all the necessary ingredients in characteristic 2 have been established previously in this paper.
\end{remark}

\appendix

\section{Two specific results about Pfister neighbors}

In this appendix, we prove two results that are needed to conclude our discussions in Remark\;\ref{5p4}.

As before, let $k$ be a field of characteristic 2. Let $\tau$ be a nonsingular  quadratic form of dimension 4 over $k$.

\newpage

\begin{prop}\label{qp2} Let $\varphi=\left<a\right>\perp [b, c] \perp \tau$, with $a,\,b,\,c\in k^*$, and suppose that $\vp$ is anisotropic.
Let $K=k(x)$ be a one-variable rational function field over $k$, and let $\psi=[f, c]\perp  \tau$, where $f=(a+bx^2)x^{-2}$.

Then $\varphi$ is a Pfister neighbor over $k$ if and only if $\psi$ is a Pfister neighbor over $K$.
\end{prop}

\begin{proof} Without loss of generality we may suppose $c=1$. Put $\theta=\left<c\right>\perp \tau=\left<1\right>\perp \tau$ and write the underlying $k$-vector space of $\theta$ as $k.s\perp U$, where $U$ denotes the vector space of $\tau$ and $\theta(s)=1$.

 Suppose that $\psi$ is a Pfister neighbor of a $3$-fold Pfister form $\pi$ over $K$. Then, $\theta_K$ is also a Pfister neighbor of $\pi$ since it is a subform of $\psi$ of dimension $5$. It follows from \cite[Prop.\;3.2 (3)]{L} that $\ind\,C_0(\theta_K)\leq 2$. Since $K/k$ is purely transcendental, we get $\ind\,C_0(\theta)\leq 2$. Again, by \cite[Prop.\;3.2 (3)]{L}, $\theta$ is a Pfister neighbor of a $3$-fold Pfister form $\rho$ defined over $k$. Consequently, $\pi\cong \rho_K$ because $\theta_K$ is a Pfister neighbor of both forms $\pi$ and $\rho_K$. Since Pfister forms are round, we get $\psi=[f, 1] \perp \tau \subset \rho_K$. Our aim is to prove that $\varphi$ is a subform of $\rho$.

 Let $W$ be the underlying $k$-vector space of $\rho$ and let
 \[
 V=\{w\in W\,|\,B_{\rho}(w,\,U)=0\}\,.
 \] The condition that $\psi$ is a subform of $\rho_K$ yields the existence of a vector $v\in V\otimes_k K$ such that $\rho(v)=f$ and $B_{\rho_K}(s,\,v)=1$.
Put $t=x^{-1}$. Note that  $f=b+ at^2\in k[t]$ and $K=k(x)=k(t)$. Applying \cite[Prop.\;17.9]{EKM08} to the anisotropic form $\rho|_V$ and the vector $s\in V$,  we can find a vector $w\in V\otimes_kk[t]$ such that $\rho_K(w)=f$ and $B_{\rho_K}(s,\,w)=1$. Since $\rho$ is anisotropic and $\deg(f)=2$, we have $w=w_0+x^{-1}w_1$ for suitable $w_0, w_1\in W$. Now it is easy to see that the following properties over $k$ hold:

\begin{itemize}
\item $B_{\rho}(w_0,\,U)=B_{\rho}(w_1,\,U)=B_{\rho}(s,\,U)=0$.
\item $B_{\rho}(s, w_0)=1$, $B_{\rho}(s, w_1)=0$,
\item $\rho(w_0)=b$, $\rho(w_1)=a$ and $B_{\rho}(w_0, w_1)=0$.
\end{itemize}
All these properties mean that the restriction of $\rho$ to the subspace $k.w_1\perp (k.s+k.w_0)\perp U$ is nothing but the form $\varphi$. Hence, $\varphi$ is a subform of $\rho$, as desired.

Conversely, suppose that $\varphi$ is a Pfister neighbor of a $3$-fold Pfister form $\rho$. Note that$$\varphi_K\cong \left<ax^{-2}\right>\perp [b, c]\perp \tau_K\cong \left<a\right>\perp [b+ax^{-2}, c] \perp \tau_K.$$In particular, $[b+ax^{-2}, c]\perp \tau_K=[f, c]\perp \tau_K$ is a Pfister neighbor of $\rho_K$. This completes the proof.
\end{proof}


\begin{prop}\label{qp3}
Let $\varphi=[a_0, a_1]\perp [b, c] \perp \tau$, with $a_0,\,a_1,\,b,\,c\in k^*$, and suppose that $\vp$ is anisotropic.
Let $L=k(t,x)$  be a two-variable rational function field over $k$, and let $\psi=[f, c]\perp  \tau$, where $f=(a_0t^2+t+a_1+bx^2)x^{-2}$.

Then $\varphi$ is a Pfister neighbor over $k$ if and only if $\psi$ is a Pfister neighbor over $L$.
\end{prop}

\begin{proof}
Suppose that $\psi$ is a Pfister neighbor over $L$. It follows from Proposition \ref{qp2} that $\varphi':=\left<a_0t^2+t+a_1\right>\perp [b, c] \perp \tau$ is a Pfister neighbor of a $3$-fold Pfister form $\pi$ over $F:=k(t)$. In particular, $[b, c] \perp \tau$ is a Pfister neighbor over $F$. Using \cite[Prop.\;3.2 (3)]{L} and arguing similarly  as in the proof of Prop.\;\ref{qp2}, we deduce that $[b, c] \perp \tau$ is a Pfister neighbor of a $3$-fold Pfister form $\rho$ over $k$. Moreover, we have $\pi\cong \rho_{F}$. We may suppose $c=1$. Hence, $\varphi'$ is a subform of $\rho_{F}$.

Let $U$ and $W$ be the underlying $k$-vector spaces of $[b, c] \perp \tau$ and $\rho$ respectively, and let $V$ be the orthogonal complement of $U$ in $W$ with respect to $B_{\rho}$. Then, as in the proof of Prop.\;\ref{qp2}, we may apply \cite[Prop.\;17.9]{EKM08} to get a vector $w\in V\otimes_kk[t]$ such that $\rho_{F}(w)=a_0t^2+t+a_1$.
Using the same arguments as in the proof of Prop.\;\ref{qp2}, there exist $w_0, w_1\in W$ such that
\[
B_{\rho}(w_i,\,U)=0\,,\; \rho(w_i)=a_i \text{ for } i=1,\,2
\]and $B_{\rho}(w_0, w_1)=1$.
These relations imply that $\varphi$ is a subform of $\rho$.

Conversely, if $\varphi$ is a Pfister neighbor, then $\varphi_{k(t)}$ is also a Pfister neighbor. Since $[a_0, a_1]_{k(t)}$ represents $a_0t^2+t+a_1$, the form $\left<a_0t^2+t+a_1\right>\perp ([b, c] \perp \tau)_{k(t)}$ is a Pfister neighbor. It follows from Prop.\;\ref{qp2} that $\psi$ is a Pfister neighbor over $L=k(t, x)$.
\end{proof}

\noindent \emph{Acknowledgements.} We thank the anonymous referee for having read the manuscript carefully and given helpful comments. In particular, Remark\;\ref{5p7} is kindly suggested by the referee. Yong Hu  is supported by grants from the National Natural Science Foundation of China (No. 11801260) and the Guangdong Basic and Applied Basic Research Foundation (No. 2021A1515010396). Peng Sun is supported by the National Key R\&D Program of China Grant No.\,2021YFA1001400 and the Fundamental Research Funds for the Central Universities.

\addcontentsline{toc}{section}{\textbf{References}}

\bibliographystyle{alpha}

\bibliography{ChowGroup}

\

Contact information of the authors:

\

Yong HU

\medskip

Department of Mathematics

Southern University of Science and Technology


Shenzhen 518055, China


Email: huy@sustech.edu.cn

\

Ahmed LAGHRIBI

\medskip

Univ. Artois, UR 2462, Laboratoire de Math\'ematiques de Lens (LML)

F-62300 Lens, France

Email: ahmed.laghribi@univ-artois.fr

\

Peng SUN

\medskip

School of Mathematics

Hunan University

Changsha 410082, China

Email: sunpeng@hnu.edu.cn

\clearpage \thispagestyle{empty}

\end{document}